\documentclass[11pt,final]{ndrpt}

\drafthead{DRAFT: \today}
\draftfoot{}

\usepackage{amssymb}
\usepackage{eepic}
\usepackage[section]{ndthms}
\usepackage[tiny,notcite,notref]{showkeys-nd}

\let\isom=\cong
\let\cong=\equiv

\let\equiv=\sim
\def\isomto{\mathrel{\hbox{$\to$\kern-.85em\raise1ex\hbox{{$\scriptstyle \isom$}}}}\;}
\def\ndiv{{\not \kern -.05em |\ }}

\let\comp=\circ

\expandafter\ifx\csname Bbb\endcsname\relax

\else

\fi

\def\arrow{\mathbin{\rightarrow}}

\def\eqdef{=_{\mathrm{df}}}

\def\otype{o}

\title{Transfinite iteration functionals and ordinal arithmetic}
\author{N. Danner}
\date{23 August 1999 ($Revision: 1.1 $)}

\def\MTpStr{\mathrm{Tp}_{\mathrm{mon}}} 
\def\TpStr{\mathrm{Tp}}			
\def\MonOmega{\Omega^{\mathrm{mon}}}

\def\hpeq{\mathbin{=^{\mathrm{hp}}}}
\def\hple{\mathbin{\leq^{\mathrm{hp}}}}

\def\mIter{I_}				

\renewcommand{\qedsymbol}{$\blacksquare$}

\newcommand{\seq}[2][\relax]{
  \ifx#1\relax
    \langle#2\rangle
  \else
    \ifx#1\left
      \left\langle#2\right\rangle
    \else
      \csname #1l\endcsname\langle#2\csname #1r\endcsname\rangle
    \fi
  \fi
}

\newcommand{\seqidx}[3][\relax]{\seq[#1]{#2}_{{#3}}}
\newcommand{\set}[2][\relax]{
  \ifx#1\relax
    \{#2\}
  \else
    \ifx#1\left
      \left\{#2\right\}
    \else
      \csname #1l\endcsname\{#2\csname #1r\endcsname\}
    \fi
  \fi
}
\newcommand{\setst}[3][\relax]{\set[#1]{#2\mid#3}}
\newcommand{\setidx}[3][\relax]{\set[#1]{#2}_{{#3}}}

\begin{document}

\begin{abstract}
Although transfinite iteration functionals have been used in the
past to construct ever-larger initial segments of the ordinals
(\cite{Feferman:Hered-replete},\cite{Aczel:Fnals-transf-type}),
there appears to be little investigation into the nature of the
functionals themselves.  In this note, we investigate the relationship
between (countable) transfinite iteration and ordinal arithmetic.
While there is a nice connection between finite iteration and addition,
multiplication, and exponentiation, we show that this it is
lost when passing to the transfinite and investigate a new equivalence
relation on ordinal functionals with respect to which we restore it.
\end{abstract}

\maketitle

\section{Introduction}

The use of functionals of higher type for defining ever-increasing
initial segments of (countable) ordinals is not a new idea---Feferman
uses a notion of transfinite iteration functionals of finite
type in \cite{Feferman:Hered-replete}
and Aczel extends this work to transfinite type
in \cite{Aczel:Fnals-transf-type}.  However,
in none of this research does there appear to be an analysis of the
iteration functionals themselves.  Specifically, we wish to understand
more completely the relationship between iteration and (ordinal)
arithmetic.  Furthermore, our original motivation for this investigation
was an interest in understanding definability of ordinals when the
tools for functional definition are restricted (this is the subject
of the author's Ph.D.\ thesis~\cite{Danner:Thesis}).  A natural way to
implement such restrictions 
is to use some version of a typed $\lambda$-calculus;
doing so necessitates that our work must take place in a structure
that can be used as a model for at least the simply typed $\lambda$-calculus.

We first consider finite iteration to determine just what such an
analysis should yield.  When we identify numbers with iterators
(for example, by representing numbers as Church numerals in the
$\lambda$-calculus),
we make explicit the view that the functional equivalent of counting
is iterated function application.  Considering counting to be the
basic operation in the universe of numbers, we are led to ask what
the numeric analogue of the basic operation of functionality is under
this equivalence.  That basic operation is, of course, application.
In other words, to what
does the interaction between iteration and application correspond in
the universe of numbers?  The most elementary interaction consists of
iterating a function, say $m$ times, then iterating it again, say $n$
times.  The result, of course, is the same as iterating the function
$m+n$ times.  In other words, application at the
object level corresponds to addition:  writing
$\mIter m^\sigma$ for the type-$\sigma$ $m$-fold iteration functional,
we have
$\mIter n^\sigma f(\mIter m^\sigma fx) = \mIter{m+n}^\sigma fx$ (associating
application to the left).  
Since iteration is defined as a higher-type functional, two more
kinds of application are basic:  application at function level and
application of one iteration functional to another.  The results
of such applications are easy to establish:
$\mIter n^\sigma(\mIter m^\sigma f) =
\mIter{mn}^\sigma f$ and $\mIter n^{\sigma\arrow\sigma}(\mIter m^\sigma) =
\mIter{m^n}^\sigma$.
Thus the fundamental operation of functionality
translates back to the universe of numbers as fundamental
operations of arithmetic:  addition, multiplication, and
exponentiation.  By viewing countable ordinals as being obtained
by transfinitely counting, the identification of numbers with
iterated function application extends to identifying countable ordinals
with transfinite iteration.
As such, we expect to see the correspondence between
application and ordinal addition, multiplication, and exponentiation extend to
transfinite iteration, and the purpose of this note is to investigate
in what way it does so.

As already
mentioned, such functionals have
been used in the past, most notably in connection with defining
ever-larger initial segments of the constructive ordinals.  
However, such work has mostly
focused on the definable ordinals, rather than the iteration
functionals themselves.  Moreover, the intuitive definition of
$\mIter\omega fx$ as $\lim_{n\to\omega}\set{\mIter nfx}$ is not
well-defined for all arguments $f$.  In order to compensate for this,
authors have usually taken the
$\omega$-iterate of a function $f$ at $x$ to be
$\sup_{n\in\omega}\set{\mIter nfx}$.  Although these definitions are
equivalent
for the functions used in practice to define ordinals (which are increasing), 
the supremum definition results in anomalies
when the focus is on iteration of arbitrary functions.  For example, if
$f(x) = 0$ for all $x$, then 
$\sup_{n\in\omega}\set{\mIter nf1} = 1$, whereas we would expect
$\mIter\omega f1$ to be $0$.

Here, we define $\alpha$-iterator functionals $\mIter\alpha^\rho$ for
each finite type $\rho$ by using 
the $\limsup$ operator,
thus staying as close as possible to the ideal of limit behavior while
maintaining totality of the functionals.  We show in
Section~\ref{sec:hm-functionals}
that if we restrict attention to monotone functions, iteration
corresponds exactly to ordinal arithmetic, as we insist
(Thm.~\ref{clm:comp-of-miter-mon}).  However,
these results cannot be extended to non-monotone functions or higher
type levels---for example,
we define a type-2 monotone functional $\Phi$
such that in general, $\mIter\gamma(\mIter\alpha\Phi)\not=
\mIter{\alpha\gamma}\Phi$.  The crux of the difficulty (which also arises when
supremum is used in the definition of $\mIter\omega$) 
is that unless the limit of a sequence exists one cannot control the
behavior of subsequences.  We resolve this in Section~\ref{sec:hp-functionals}
by introducing a new equivalence relation $\hpeq$ on
ordinal functionals
which allows us to focus our attention on
arguments for which the appropriate limits do exist (although, as
mentioned above, we do not eliminate such arguments from consideration
altogether).  We then establish the desired correspondence relative
to $\hpeq$ for all functionals at all type levels
(Thm.~\ref{thm:comp-of-miter}).  As $\hpeq$ is just
equality on the ordinals themselves, we can make use of the
correspondence to define larger ordinals through application of
iteration functionals.

\section{Preliminaries}
\label{sec:Preliminaries}

We will work in two finite type structures over $\Omega$,
where $\Omega$ is the first uncountable ordinal.  We define the full
type structure
$\TpStr(\Omega) = \setidx{\Omega_\sigma}{\sigma}$ and the hereditarily
monotone type structure
$\MTpStr(\Omega) = \setidx{\MonOmega_\sigma}{\sigma}^{\vphantom{\mathrm{mon}}}$
as follows.
$\Omega_\otype^{\vphantom{\mathrm{mon}}} = 
\MonOmega_\otype = \Omega$, and the order in both
cases is the usual order on the ordinals.  If $\Omega_\sigma$ and
$\Omega_\tau$ have been defined, then
\begin{align*}
\Omega_{\sigma\arrow\tau} &=
 \setst{f}{f:\Omega_\sigma\to\Omega_\tau} \\
\MonOmega_{\sigma\arrow\tau} &=
 \setst{f}{f:\MonOmega_\sigma\to\MonOmega_\tau\text{ is monotone}}
\end{align*}
where we say that $f$ is \emph{monotone} provided that $f(x)\leq f(y)$
whenever $x\leq y$.
The order is defined pointwise in both cases:  $f\leq g$ if for all
$x\in\Omega_\sigma^{\vphantom{\mathrm{mon}}}$ ($x\in\MonOmega_\sigma$),
$f(x)\leq g(x)$.

The pointwise definition of the order on $\Omega_{\sigma\arrow\tau}$
yields a pointwise characterization of supremums and infimums over
an arbitrary index set $I$:
$$
  \Bigl(\sup_{i\in I}\set{f_i}\Bigr)(x) = \sup_{i\in I}\set[big]{f_i(x)}
  \qquad
  \Bigl(\inf_{i\in I}\set{f_i}\Bigr)(x) = \inf_{i\in I}\set[big]{f_i(x)}
$$

\begin{prop}
For each type $\sigma$, if $X\subseteq \Omega_\sigma$, then $\inf X$
exists; moreover, if $X$ is countable, then $\sup X$ exists.
\end{prop}
\begin{proof}
Both claims are proved by induction on $\sigma$.  The existence of
$\sup X$ in the base case follows from the regularity of $\Omega$ and
the induction step is trivial.
\end{proof}

The following definitions of $\limsup$, $\liminf$ and limit are taken
from Birkhoff~\cite[\S X.9]{Birkhoff:Lattice-Theory-3}, but we have
restricted attention to the case in which the nets are based on
countable ordinals.

\begin{defn}
For each type $\sigma$ and countable ordinal $\zeta$, if
$\setidx{x_\xi}{\xi<\zeta}\subseteq\Omega_\sigma$, then
$$
  \limsup_{\xi\to\zeta}\set{x_\xi}\eqdef\mathop{\vphantom{\sup}\inf}\limits_{\gamma<\zeta}\set[Big]{\sup_{\gamma\leq\xi<\zeta}\set{x_\xi}}
  \qquad
  \liminf_{\xi\to\zeta}\set{x_\xi}\eqdef\sup_{\gamma<\zeta}\set[Big]{\mathop{\vphantom{\sup}\inf}\limits_{\gamma\leq\xi<\zeta}\set{x_\xi}}
$$
If there is $x\in\Omega_\sigma$ such that $\limsup_{\xi\to\zeta}\set{x_\xi} =
x = \liminf_{\xi\to\zeta}\set{x_\xi}$, then we say that
$\lim_{\xi\to\zeta}\set{x_\xi}$ exists and is equal to $x$.
\end{defn}

Both Aczel \cite{Aczel:Fnals-transf-type} and the author
\cite{Danner:Thesis} have defined transfinite type structures.
In both cases, limit level function spaces are defined as a product
over the function spaces of lower level (in fact, Aczel defines successor
levels in the same way for the sake of uniformity).  By extending the order to
such spaces coordinatewise, it is not difficult to extend the results in
this paper to such type structures.

\begin{defn}
For each type $\sigma$ and countable ordinal $\alpha$, the 
$\alpha$-iteration functional of type $\sigma$ is the
functional
$\mIter\alpha^\sigma:\Omega_{\sigma\arrow\sigma}\to\Omega_{\sigma\arrow\sigma}$
defined by
$$
  \mIter 0^\sigma fx = x
  \qquad
  \mIter{\alpha+1}^\sigma fx = f(\mIter\alpha^\sigma fx)
  \qquad
  \mIter\mu^\sigma fx = \limsup_{\xi\to\mu}\{\mIter\xi^\sigma fx\}
$$
where application associates to the left.
We usually drop the type subscript whenever it is clear from 
context or irrelevant.
\end{defn}

Note that by using the pointwise characterization of $\sup$ and $\inf$, we
can push arguments in and out of $\limsup$:
$\limsup_{\xi\to\zeta}\set{f_\xi x} = (\limsup_{\xi\to\zeta}\set{f_\xi})(x)$.
Applying this to the definition of $\mIter\mu$ for limit $\mu$, we have
$\mIter\mu fx = \limsup_{\xi\to\mu}\set{\mIter\xi fx} =
(\limsup_{\xi\to\mu}\set{\mIter\xi})fx$, so $\mIter\mu =
\limsup_{\xi\to\mu}\set{\mIter\xi}$.

We give two counterexamples to show that the correspondence between
transfinite iteration and ordinal arithmetic need not hold.  Let
$f$ be any ordinal function such that $f(2x) = 1$ and $f(2x+1) = 0$
when $x<\omega$.  Then
if $g\eqdef \mIter2 f$, we have $g(2x) = 0$ and
$g(2x+1) = 1$ for all $x<\omega$,
so $\mIter\omega(\mIter 2f)(0) = \mIter\omega g0 = 0$.  On the other
hand, $\mIter{2\omega} f0 = \mIter\omega f0 = 1$.  Of course, $f$ is
a rather poorly-behaved function, and one might hope that this difficulty
would not arise for functions that are somehow well-behaved.  
For example, Aczel~\cite{Aczel:Fnals-transf-type}
restricts attention to hereditarily inflationary functionals.  
This is not an ideal resolution for us
for two reasons:  it requires a ``pure'' type structure (i.e., functionals
always have the same domain and range) so that it makes sense to compare
input and output, and such functionals do not yield a model in which
the $\lambda$-calculus can be directly interpreted 
(since, e.g., constant functionals are $\lambda$-definable but
not inflationary).

We give another example of the failure of
application to correspond to arithmetic, this time using only
monotone functionals.  In particular, we cannot equate
``well-behaved'' with monotonicity.  The
type-2 functional to be iterated interchanges two 
functions.  In this case, the double iterate will be
the identity on either of those functions, and so the $\omega$-iterate
of the double iterate will also be the identity on either of the functions.
However, the $\omega$-iterate of the functional itself cannot be the
identity, because it is alternating between the two.  
For the two functions,
define $f_0(\alpha) = \alpha$, $f_1(\alpha) = 2$ for all $\alpha$.  Set
$g\eqdef\max\{f_0,f_1\}$ and $h\eqdef\min\{f_0,f_1\}$, and
define
$\Phi$ by
$$
  \Phi(f) =
  \begin{cases}
    h,&f\leq h \\
    f_{1-i},&f\leq f_i,f\nleq h \\
    g,&\text{otherwise}
  \end{cases}
$$
Verifying that $\Phi$ is monotone is straightforward, though tedious;
the following picture of the action of $\Phi$ should suffice:
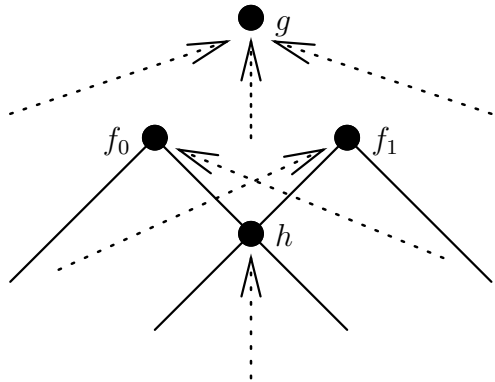
\begin{center}
\setlength{\unitlength}{0.00083333in}
\begingroup\makeatletter\ifx\SetFigFont\undefined
\def\x#1#2#3#4#5#6#7\relax{\def\x{#1#2#3#4#5#6}}%
\expandafter\x\fmtname xxxxxx\relax \def\y{splain}%
\ifx\x\y   
\gdef\SetFigFont#1#2#3{%
  \ifnum #1<17\tiny\else \ifnum #1<20\small\else
  \ifnum #1<24\normalsize\else \ifnum #1<29\large\else
  \ifnum #1<34\Large\else \ifnum #1<41\LARGE\else
     \huge\fi\fi\fi\fi\fi\fi
  \csname #3\endcsname}%
\else
\gdef\SetFigFont#1#2#3{\begingroup
  \count@#1\relax \ifnum 25<\count@\count@25\fi
  \def\x{\endgroup\@setsize\SetFigFont{#2pt}}%
  \expandafter\x
    \csname \romannumeral\the\count@ pt\expandafter\endcsname
    \csname @\romannumeral\the\count@ pt\endcsname
  \csname #3\endcsname}%
\fi
\fi\endgroup
\begin{picture}(3044,2369)(0,-10)
\thicklines
\put(1522,2271){\blacken\ellipse{150}{150}}
\put(1522,2271){\ellipse{150}{150}}
\dottedline{90}(3022,1671)(1672,2121)
\path(1918.658,2102.026)(1672.000,2121.000)(1880.710,1988.184)
\dottedline{90}(22,1671)(1372,2121)
\path(1163.290,1988.184)(1372.000,2121.000)(1125.342,2102.026)
\dottedline{67}(1522,1521)(1522,2121)
\path(1582.000,1881.000)(1522.000,2121.000)(1462.000,1881.000)
\put(1672,2196){\makebox(0,0)[lb]{\smash{{{\SetFigFont{12}{14.4}{rm}$g$}}}}}
\put(922,1522){\blacken\ellipse{150}{150}}
\put(922,1522){\ellipse{150}{150}}
\put(1522,922){\blacken\ellipse{150}{150}}
\put(1522,922){\ellipse{150}{150}}
\put(2122,1522){\blacken\ellipse{150}{150}}
\put(2122,1522){\ellipse{150}{150}}
\path(22,622)(922,1522)(2122,322)
\path(3022,622)(2122,1522)(922,322)
\dottedline{90}(322,697)(1972,1447)
\path(1778.340,1293.065)(1972.000,1447.000)(1728.684,1402.309)
\dottedline{90}(2722,772)(1072,1447)
\path(1316.849,1411.661)(1072.000,1447.000)(1271.413,1300.595)
\dottedline{90}(1522,22)(1522,772)
\path(1582.000,532.000)(1522.000,772.000)(1462.000,532.000)
\put(2272,1447){\makebox(0,0)[lb]{\smash{{{\SetFigFont{12}{14.4}{rm}$f_1$}}}}}
\put(1672,847){\makebox(0,0)[lb]{\smash{{{\SetFigFont{12}{14.4}{rm}$h$}}}}}
\put(772,1447){\makebox(0,0)[rb]{\smash{{{\SetFigFont{12}{14.4}{rm}$f_0$}}}}}
\end{picture}
\end{center}
We want to compare $\mIter\omega(\mIter2\Phi)f_0$ and
$\mIter\omega\Phi f_0$.  For the former, set
$\Psi\eqdef\mIter2\Phi$; then
$$
  \Psi(f) =
  \begin{cases}
    h,&f\leq h \\
    f_i,&f\leq f_i,f\nleq h \\
    g,&\text{otherwise}
  \end{cases}
$$
In particular, $\Psi(f_i) = f_i$ for $i=0,1$,
so $\mIter\omega(\mIter2\Phi)f_0 = \mIter\omega\Psi f_0 = f_0$.
On the other hand, a direct computation shows that
$\mIter\omega\Phi f_0 = g$, and therefore
$\mIter\omega(\mIter2\Phi)f_0\not=\mIter\omega\Phi f_0$.

What drives this example is the fact that $f_0$ and
$\Phi(f_0)$ are not comparable---as a result,
the sequence of iterates $\seqidx{\Phi^n(f)}{n}$ does not have a
limit, and therefore subsequences may have different limiting behavior 
than the sequence.  We begin to
repair the damage by analyzing
iteration of monotone functionals which map
each input to a comparable output.  In this case the iterates form
either non-decreasing or non-increasing
sequences, and as a result subsequences will behave well.
Unfortunately, the
comparability requirement is too restrictive, because it is only guaranteed
to hold when the order on the domain is total.  Thus, it
prevents us from establishing the connection between application
of the iteration functionals and arithmetic at 
higher type.  To push upwards,
we develop the notion of hereditarily positive
equality, with respect to which the correspondence is exact at all types.

\section{Hereditarily Monotone Functionals}
\label{sec:hm-functionals}

We partially investigated iteration functionals in $\MTpStr(\Omega)$
in~\cite{Danner:Ords-ord-fns}; the results here significantly
extend this earlier work.

\begin{lemma}
\label{clm:hm-closed-under-limsup}
For each type $\sigma$ and countable ordinal $\zeta$,
if $f_\xi\in\MonOmega_\sigma$
for all $\xi<\zeta$, then
$\limsup_{\xi\to\zeta}\set{f_\xi}\in\MonOmega_\sigma$.
\end{lemma}
\begin{proof}
The Lemma is proved
by induction on $\sigma$.  
This is trivial if $\MonOmega_\sigma = \MonOmega$.  Otherwise,
suppose that $\sigma = \rho\arrow\tau$.  If $x\in\MonOmega_\rho$, then
$(\limsup\set{f_\xi})(x) = \limsup\set{f_\xi x}\in\MonOmega_\tau$ by induction,
because each $f_\xi x\in\MonOmega_\tau$.  Furthermore, 
if $x\leq x'$ are elements of
$\MonOmega_\rho$, then $(\limsup\set{f_\xi})(x) =
\limsup\set{f_\xi x} \leq\limsup\set{f_\xi x'} = 
(\limsup\set{f_\xi})(x)$, with the inequality holding because $f_\xi$
is hereditarily monotone and $x\leq x'$, so $f_\xi x\leq f_\xi x'$ for 
all $\xi$.
\end{proof}

\begin{prop}
For each type $\sigma$ and countable ordinal $\alpha$, 
$\mIter\alpha^\sigma$ is hereditarily monotone.
\end{prop}
\begin{proof}
The Proposition is proved
by induction on $\alpha$ for all $\sigma$.  
If $\alpha = 0$, then $\mIter\alpha$ is just
the functional that is constantly the identity on $\MonOmega_\sigma$, 
which is easily seen to
be hereditarily monotone.  Suppose that $\alpha = \gamma+1$.  First we must
verify that $\mIter\alpha$ maps $\MonOmega_{\sigma\arrow\sigma}$ to itself.
Suppose that  $f\in\MonOmega_{\sigma\arrow\sigma}$ and
$x\in\MonOmega_\sigma$.  Then since $\mIter\gamma fx\in\MonOmega_\sigma$
by the induction hypothesis and $f$ is hereditarily
monotone by assumption, $\mIter\alpha fx = f(\mIter\gamma fx)$ is 
hereditarily monotone,
and so $\mIter\alpha f\in\MonOmega_{\sigma\arrow\sigma}$.
We must also verify that if $x$ and $x'$ are hereditarily monotone,
$x\leq x'$, then $\mIter\alpha fx\leq\mIter\alpha fx'$, which is just as
easy to do.  Second, we must verify the monotonicity of
$\mIter\alpha$:  if $f$, 
$f'\in\MonOmega_{\sigma\arrow\sigma}$ are such that $f\leq f'$,
then $\mIter\alpha f\leq\mIter\alpha f'$.
Fix any $x\in\MonOmega_\sigma$.  Then $\mIter\alpha fx =
f(\mIter\gamma fx) \leq f(\mIter\gamma f' x) \leq f'(\mIter\gamma f'x) =
\mIter\alpha f'x$; the first inequality follows from the fact that
$\mIter\gamma f\leq \mIter\gamma f'$ (induction) and the second from the
fact that $f\leq f'$.  This takes care of the successor case.  If $\alpha$
is a limit, then $\mIter\alpha = \limsup_{\xi\to\alpha}\set{\mIter\xi}$
by definition; but this $\limsup$ is hereditarily monotone by induction
and Lemma~\ref{clm:hm-closed-under-limsup}.
\end{proof}

\begin{defn}
We say that $\setidx{f_\xi}{\xi<\zeta}\subseteq\MonOmega_\sigma$ is
\emph{non-decreasing} 
(\emph{non-increasing}) if whenever
$\alpha<\gamma<\zeta$, $f_\alpha\leq f_\gamma$
($f_\alpha\geq f_\gamma$, respectively).  We use the same terminology
when $\setidx{f_\xi}{\xi<\zeta}\subseteq\Omega_\sigma$.
\end{defn}

\begin{lemma}
\label{clm:lim-of-seq}
Fix any type $\sigma$ and countable ordinal $\zeta$, and
let $\setidx{f_\xi}{\xi<\zeta}\subseteq\MonOmega_\sigma$.
\begin{enumerate}
\item \label{item:lim-of-subseq}
If the limit of a sequence exists,
then it is the limit of any subsequence:
if $\lim\set{f_\xi}$ exists, $q:\zeta'\to\zeta$ is non-decreasing, and
$\lim_{\alpha\to\zeta'}q(\alpha) = \zeta$, then
$\lim\setidx{f_{q(\xi)}}{\xi<\zeta'}$ exists and is equal to
$\lim\setidx{f_\xi}{\xi<\zeta}$.
\item \label{item:lim-non-decr}
If $\set{f_\xi}$ is non-decreasing, then $\lim\set{f_\xi}$ exists and
is equal to $\sup_{\xi<\zeta}\set{f_\xi}$.
\item \label{item:lim-non-incr}
If $\set{f_\xi}$ is non-increasing, then $\lim\set{f_\xi}$ exists and
is equal to $\inf_{\xi<\zeta}\set{f_\xi}$.
\end{enumerate}
\end{lemma}
\begin{proof}
(\ref{item:lim-of-subseq}) In general, the $\liminf$ of a subsequence is
always greater than or equal to the $\liminf$ of the sequence, and
vice-versa for $\limsup$, so if $f = \lim_{\xi\to\zeta}\set{f_\xi}$, then
$\liminf\set{f_{q(\xi)}}\geq f \geq \limsup\set{f_{q(\xi)}}$.  But for
any sequence $\setidx{g_\mu}{\mu<\theta}$, 
$\liminf\set{g_\mu}\leq\limsup\set{g_\mu}$, so this implies that
$\liminf{f_{q(\xi)}} = f = \limsup_{f_{q(\xi)}}$.

(\ref{item:lim-non-decr}) Since $\set{f_\xi}$ is non-decreasing,
$$
  \liminf_{\xi\to\zeta}\set{f_\xi} = 
  \sup_{\mu<\zeta}\set[big]{\inf_{\mu\leq\xi<\zeta}\set{f_\xi}} = 
  \sup_{\mu<\zeta}\set{f_\mu}
$$
and
$$
  \limsup_{\xi\to\zeta}\set{f_\xi} =
  \inf_{\mu<\zeta}\set[big]{\sup_{\mu\leq\xi<\zeta}\set{f_\xi}} =
  \sup_{0\leq\xi<\zeta}\set{f_\xi} =
  \liminf_{\xi\to\zeta}\set{f_\xi}.
$$

(\ref{item:lim-non-incr}) is similar to (\ref{item:lim-non-decr}).
\end{proof}

\begin{lemma}
\label{clm:mon-on-sup-inf}
Fix any types $\sigma$ and $\tau$, $X\subseteq\MonOmega_\sigma$, and
let $f:\MonOmega_\sigma\to\MonOmega_\tau$ be monotone.
If $\sup X$ exists, then $f(\sup X)\geq\sup\setst{f(x)}{x\in X}$, and
if $\inf X$ exists, then $f(\inf X)\leq\inf\setst{f(x)}{x\in X}$.
\end{lemma}
\begin{proof}
Both claims have similar proofs, so we just do the first.  If $x\in X$,
then $x\leq\sup X$, so by monotonicity of $f$, $f(x)\leq f(\sup X)$.  Since
$x$ was chosen arbitrarily, $\sup\setst{f(x)}{x\in X}\leq f(\sup X)$.
\end{proof}

\begin{lemma}
\label{clm:iterates-monotone}
For each type $\sigma$, countable ordinal $\zeta$,
$f\in\MonOmega_{\sigma\arrow\sigma}$, and $x\in\MonOmega_\sigma$:
\begin{enumerate}
\item If $\mIter{\zeta+1}fx\geq \mIter\zeta fx$, then for all
$\gamma>\alpha\geq\zeta$, $\mIter\gamma fx\geq\mIter\alpha fx$.
\item If $\mIter{\zeta+1}fx\leq \mIter\zeta fx$, then for all
$\gamma>\alpha\geq\zeta$, $\mIter\gamma fx\leq\mIter\alpha fx$.
\end{enumerate}
\end{lemma}
\begin{proof}
Each clause is proved by a similar
induction on $\gamma$; we do just the first.  
Throughout the proof, we make silent use of 
Lemma~\ref{clm:lim-of-seq}(\ref{item:lim-of-subseq}) to identify the limit
of a sequence with the limit of a tail of that sequence, provided the
former exists.
If $\gamma=0$, the
claim is vacuous.  Suppose that $\gamma = \delta+1$; by induction, it
suffices to show that $\mIter\delta fx\leq\mIter\gamma fx$, and we do this
by induction on $\delta$.  If $\delta =\alpha$, 
then this is just the hypothesis
that $\mIter\zeta fx\leq \mIter{\zeta+1}fx$.  
The successor case is straightforward.  Suppose that
$\delta$ is a limit.  By the main induction hypothesis, 
$\setidx{\mIter\xi fx}{\zeta\leq\xi<\delta}$ is a non-decreasing sequence, so
$\mIter\delta fx = \lim_{\xi\to\delta}\set{\mIter\xi fx} = 
\sup_{\zeta\leq\xi<\delta}\set{\mIter\xi fx}$ by
Lemma~\ref{clm:lim-of-seq}(\ref{item:lim-non-decr}).  Now applying
Lemma~\ref{clm:mon-on-sup-inf}, 
$$
\mIter\gamma fx = f(\mIter\delta fx) = 
f\bigl(\sup_{\zeta\leq\xi<\delta}\set{\mIter\xi fx}\bigr) \geq \\
\sup_{\zeta\leq\xi<\delta}\set{f(\mIter\xi fx)} =
\sup_{\zeta\leq\xi<\delta}\set{\mIter{\xi+1} fx}.
$$
This last sequence is a subsequence of
$\setidx{\mIter\xi fx}{\zeta\leq\xi<\delta}$, 
so it is non-decreasing, and therefore by 
Lemma~\ref{clm:lim-of-seq}(\ref{item:lim-non-decr})
its supremum is a limit, and by
Lemma~\ref{clm:lim-of-seq}(\ref{item:lim-of-subseq})
the limit is the same as that of the original sequence:
$\sup_{\zeta\leq\xi<\delta}\set{\mIter{\xi+1}fx}=\lim\set{\mIter{\xi+1}fx} =
\lim\set{\mIter\xi fx} = \mIter\delta fx$.
So $\mIter\gamma fx\geq \mIter\delta fx$.  This completes the induction step
for successor $\gamma$.  Finally, suppose that $\gamma$
is a limit.  Then by induction $\setidx{\mIter\xi fx}{\xi<\gamma}$ is
non-decreasing, so for any $\alpha<\gamma$,
$\mIter\alpha fx \leq\sup_{\xi<\gamma}\set{\mIter\xi fx} =
\mIter\gamma fx$.
\end{proof}

\begin{prop}
\label{clm:miter-limit-mon}
If $f:\MonOmega\to\MonOmega$ is a monotone function
and $\alpha$ is a countable limit ordinal, then for all $\beta$,
$\mIter\alpha f\beta = \lim_{\xi\to\alpha}\mIter\xi f\beta$.
\end{prop}
\begin{proof}
This follows from Lemmas~\ref{clm:lim-of-seq} and \ref{clm:iterates-monotone}
(taking $\zeta = 0$),
because the order on $\MonOmega$ is total.
\end{proof}

We can now establish the connection between arithmetic of
ordinals and application of iteration functionals at base type:

\begin{thm}[Iteration Functionals in $\MTpStr(\MonOmega)$]
\label{clm:comp-of-miter-mon}
Suppose $f:\MonOmega\to\MonOmega$ is a monotone function.
Then for any $\alpha$ and $\gamma$:
\begin{enumerate}
\item \label{item:mon-add}
$\mIter\alpha^\Omega f\comp\mIter\gamma^\Omega f = 
\mIter{\gamma+\alpha}^\Omega f$.
\item \label{item:mon-mult}
$\mIter\alpha^\Omega(\mIter\gamma^\Omega f) = 
\mIter{\gamma\alpha}^\Omega f$.
\item \label{item:mon-exp}
$\mIter\alpha^{\Omega\arrow\Omega}(\mIter\gamma^\Omega)f =
\mIter{\gamma^\alpha}^\Omega f$.
\end{enumerate}
\end{thm}
\begin{proof}
All three clauses
are proved by induction on $\alpha$; we do (\ref{item:mon-mult}) as an
example.  Fix any ordinal $\beta$.  If 
$\alpha = 0$, then $\mIter\alpha(\mIter\gamma f)\beta = \beta =
\mIter{\alpha\gamma}f\beta$.  

If $\alpha = \delta+1$, then
$\mIter\alpha(\mIter\gamma f)\beta
= \mIter\gamma f(\mIter\delta(\mIter\gamma f)\beta)
= \mIter\gamma f(\mIter{\gamma\delta}f\beta)
= \mIter{\gamma\delta+\gamma}f\beta
= \mIter{\gamma\alpha}f\beta$, where the second equality follows from
the induction hypothesis and the third from part (\ref{item:mon-add}).

Suppose that $\alpha$ is a limit.
By Prop.~\ref{clm:miter-limit-mon}, $\mIter{\gamma\alpha}f\beta =
\lim_{\xi\to\gamma\alpha}\set{\mIter\xi f\beta}$.  Since
$\setidx{\mIter{\gamma\xi}f\beta}{\xi<\alpha}$ is a subsequence of
$\setidx{\mIter\xi f\beta}{\xi<\gamma\alpha}$ and the limit of the
latter sequence exists,
\begin{align*}
\mIter{\gamma\alpha}f\beta
 &= \lim_{\xi\to\gamma\alpha}\set{\mIter\xi f\beta} 
   & \qquad & \text{(Prop.~\ref{clm:miter-limit-mon})} \\
 &= \lim_{\xi\to\alpha}\set{\mIter{\gamma\xi}f\beta} 
   & \qquad & \text{(Lemma~\ref{clm:lim-of-seq}(\ref{item:lim-of-subseq}))} \\
 &= \limsup_{\xi\to\alpha}\set{\mIter{\gamma\xi}f\beta} 
   & \qquad & \text{(Definition of $\lim$)} \\
 &= \limsup_{\xi\to\alpha}\set{\mIter\xi(\mIter\gamma f)\beta}
    & \qquad & \text{(Induction Hypothesis)} \\
 &= \Bigl(\limsup_{\xi\to\alpha}\set{\mIter\xi}\Bigr)(\mIter\gamma f)(\beta) 
   & \qquad & \text{(Definition of $\limsup$)} \\
 &= \mIter\alpha(\mIter\gamma f)\beta
    & \qquad & \text{(Definition of $\mIter\alpha$)}
\end{align*}
completing the proof.
\end{proof}

We show by example that the hypothesis of Lemma~\ref{clm:iterates-monotone}
need not be satisfied at higher type.
It suffices to find a monotone function $f$ such that 
$\mIter\gamma^{\Omega\arrow\Omega} f$ is not 
comparable with $f$ for some $\gamma$.
Consider the function $f$ defined by:
$$
  f(\xi) =
  \begin{cases}
  \xi+1,&\xi<\omega \\
  \omega,&\xi = \omega,\xi = \omega+1 \\
  \omega+1,&\xi>\omega+1
  \end{cases}
$$
Then $f$ is monotone, but $\mIter\omega f$ is the function that is constantly
$\omega$, so $\mIter\omega f$ is not comparable with $f$.  
We also recall that we showed with the functional $\Phi$ in
the previous section that
we cannot extend Thm.~\ref{clm:comp-of-miter-mon}\footnote{Actually, 
it is possible to extend
part (\ref{item:mon-add}) by using the fact that for any
$\mu<\alpha$, $\limsup_{\xi\to\alpha}\set{f_\xi} =
\limsup_{\mu<\xi\to\alpha}\set{f_\xi}$.} to the type
$\MonOmega\arrow\MonOmega$.

\section{Hereditarily Positive Functionals}
\label{sec:hp-functionals}

In order
to establish the desired correspondence between application of iteration
functionals and arithmetic at higher type, we introduce a new
notion: hereditarily positive equality.  However, the
result that we
prove (Thm.~\ref{thm:comp-of-miter})
is technically weaker than Thm.~\ref{clm:comp-of-miter-mon}
and cannot be used to derive the latter.  Nonetheless, as the new
equivalence relation is just equality on the ordinals, it is sufficient
for defining them.
In this section, we work in the full type structure
$\TpStr(\Omega)$.

\begin{defn}
The \emph{hereditarily positive} 
(\emph{h.p.}) functionals and
the order $\hple$
are defined
simultaneously by induction on type as follows:
\begin{itemize}
\item Any element of $\Omega$ or $\Omega_{\rho\arrow\tau}$,
$\rho\not=\tau$, is h.p.; $\hple$ in either case is just $\leq$.
\item If $f\in\Omega_{\tau\arrow\tau}$, then $f$ is h.p.\ provided:
\begin{itemize}
\item If $x\in\Omega_\tau$ is h.p., then $fx$ is h.p.;
\item $f$ is hereditarily inflationary\footnote{We use the phrase
``hereditarily inflationary'' instead of the more accurate but somewhat
wordier ``inflationary on h.p. arguments'', and similarly we say
``hereditarily monotone''.}:
if $x\in\Omega_\tau$ is h.p., then $x\hple fx$;
\item $f$ is hereditarily monotone:
if $x$, $x'\in\Omega_\tau$ are h.p. and $x\hple x'$, then $fx\hple fx'$.
\end{itemize}
If $f$, $f'\in\Omega_{\tau\arrow\tau}$, we say $f\hple f'$ provided that
for all h.p.\ $x\in\Omega_{\tau}$, $fx\hple f'x$.
\end{itemize}
We say that $f\hpeq g$ if $f\hple g$ and $g\hple f$.
\end{defn}

We stress that the h.p.\ functionals do \emph{not} form a new 
type structure---they are a subclass of the universe of
an existing one.  However, the order $\hple$ itself is defined on \emph{all}
functionals, even those that are not themselves hereditarily positive.
When proving facts involving the notion of hereditarily positive, we will
often use induction on type---in this situation, there are two base cases:
the type $\Omega$, and all types of the form
$\Omega_{\sigma\arrow\tau}$ with $\sigma\not=\tau$.

\begin{lemma}
\label{lem:hple-facts}
~
\begin{enumerate}
\item \label{item:refl-trans}
$\hple$ is reflexive and transitive, and therefore $\hpeq$ is an
equivalence relation.
\item \label{item:eq-implies-hpeq}
If $f\leq g$, then $f\hple g$; if $f=g$, then $f\hpeq g$; if $f$, 
$f'\in\Omega_{\sigma\arrow\tau}$, then $f\hpeq f'$ iff
for all $x\in\Omega_\sigma$, $fx\hpeq f'x$.
\item \label{item:boundedness}
If $q:\zeta\to\zeta'$ and for all $\xi<\zeta$, $f_\xi\hple
f_{q(\xi)}'$, then $\limsup_{\xi\to\zeta}\{f_\xi\}\hple
\limsup_{\xi\to\zeta}\{f_{q(\xi)}'\}$.  In particular,
if for all $\xi<\zeta$, $f_\xi\hple f'_\xi$, then 
$\limsup\{f_\xi\}\hple\limsup\{f_\xi'\}$, and if $f_\xi\hple f$ for
all $\xi$, then $\limsup_{\xi\to\zeta}\{f_\xi\}\hple f$.
\item \label{item:limsup-eq-sup}
If $f_\alpha\hple f_\gamma$ for $\alpha<\gamma < \zeta$, then
$\limsup_{\xi\to\zeta}\{f_\xi\} \hpeq \sup_{\xi<\zeta}\{f_\xi\}$.
\end{enumerate}
\end{lemma}
\begin{proof}
(1) and (2) are immediate, and (3) and (4) are proved by induction on type.
We provide details for (4).
Note that this is not a trivial claim, as it is an assertion about the
h.p.\ order, not the pointwise order.
The claim is true for the base cases because the two orders are the same.
Suppose $f_\xi:\Omega_\tau\to\Omega_\tau$ for
$\xi<\zeta$.
If $\alpha<\gamma<\zeta$ and $x\in\Omega_\tau$ is h.p.,
then since $f_\alpha\hple f_\gamma$, 
we have $f_\alpha x\hple f_\gamma x$, and hence
$\left(\limsup_{\xi\to\zeta}\{f_\xi\}\right)x = 
\limsup_{\xi\to\zeta}\{f_\xi x\} \hpeq
\sup_{\xi<\zeta}\{f_\xi x\} = \left(\sup_{\xi<\zeta}\{f_\xi\}\right)x$,
with the second equality following from the induction hypothesis.
\end{proof}

\begin{lemma}
\label{lem:limsup-hp}
Fix any type $\sigma$, countable ordinal $\zeta$, and
$\setidx{f_\xi}{\xi<\zeta}\subseteq\Omega_\sigma$.
If there is $\alpha$ such that $f_\xi$ is h.p.\ for 
all $\xi\geq\alpha$, then $\limsup_{\xi\to\zeta}\set{f_\xi}$
is h.p.
\end{lemma}
\begin{proof}
The lemma follows from the special case $\alpha=0$,
since the $\limsup$ of a sequence is the same as the $\limsup$ of
any tail of that sequence.
The proof is
by induction on $\sigma$, using Lemma~\ref{lem:hple-facts}.
The claim is trivially true in the base cases.  Suppose that
$f_\xi:\Omega_\tau\to\Omega_\tau$.
\begin{itemize}
\item If $x$ is h.p., then for all $\xi$, $f_\xi x$ is h.p., so
$\left(\limsup\{f_\xi\}\right)x = \limsup\{f_\xi x\}$ is
h.p.\ by the induction hypothesis.
\item If $x$ is h.p., then for all $\xi$ we have $x\hple f_\xi x$, so
$x \hple\limsup\{f_\xi x\} =
\left(\limsup\{f_\xi\}\right)x$.
\item If $x\hple x'$ are h.p., then for all $\xi$ we have
$f_\xi x\hple f_\xi x'$, 
so $\left(\limsup\{f_\xi\}\right)x = \limsup\{f_\xi x\}
\hple \limsup\{f_\xi x'\} = \left(\limsup\{f_\xi\}\right)x'$.
\hfill\qedsymbol
\end{itemize}
\def\qed{}
\end{proof}

As $\mIter 0$ is the functional that is constantly the identity,
it is not inflationary and hence not h.p.
However, this is the only way in which the iteration functionals are
not well-behaved:
$\mIter\alpha$ is h.p.\ for all $\alpha\geq 1$, and the functionals
$\mIter\alpha$ form a non-decreasing sequence with respect to
$\hple$.

\begin{prop}
\label{clm:miter-hp}
For each type $\sigma$ and countable 
$\alpha\geq 1$, $\mIter\alpha^\sigma$ is h.p.
\end{prop}
\begin{proof}
The proof is by induction on $\alpha$.  
If $\alpha = 1$, then $\mIter\alpha$ is
the identity function, which is easily seen to be h.p.

Suppose that $\alpha=\gamma+1$.  First we must show that if $f$
is h.p., then so is $\mIter\alpha f$,
using the fact that $I_\gamma f$ is h.p.\ by the induction hypothesis.
\begin{itemize}
\item If $x$ is h.p., then $I_\alpha fx = f(I_\gamma fx)$ is h.p.\
because $I_\gamma fx$ is h.p.\ by the induction hypothesis and $f$ maps h.p.\
functionals to h.p.\ functionals by assumption.
\item If $x$ is h.p., then $x\hple I_\gamma fx$, so 
$x\hple fx\hple f(I_\gamma fx) = I_\alpha fx$.  The first inequality follows
from the fact that $f$ is hereditarily inflationary, the
second from the fact that $f$ is hereditarily monotone.
\item If $x\hple x'$ are h.p., then $I_\alpha fx = f(I_\gamma fx) 
\hple f(I_\gamma fx') = I_\alpha fx'$.  The second inequality follows from
the facts that $I_\gamma f$ is h.p. and $f$ is hereditarily monotone.
\end{itemize}
To show that $I_\alpha$ is hereditarily inflationary, it suffices to show that
if $f$ and $x$ are h.p., then $fx\hple I_\alpha fx$, which we did above.
To show that $I_\alpha$ is hereditarily monotone, fix
$f\hple f'$ and $x\hple x'$ and note that
$I_\alpha fx = f(I_\gamma fx)\hple f(I_\gamma fx')\hple f(I_\gamma f'x')
\hple f'(I_\gamma f'x') = I_\alpha f'x'$, repeatedly using the induction
hypothesis and hereditary monotonicity of h.p.\ functionals.

If $\alpha$ is a limit, then $\mIter\alpha =
\limsup_{\xi\to\alpha}\set{\mIter\xi}$ is h.p.\ by
Lemma~\ref{lem:limsup-hp} because
$\mIter\xi$ is h.p.\ for all $1\leq\xi<\alpha$ by the inductive
hypothesis.
\end{proof}

\begin{prop}
\label{prop:miter-incr}
For all countable $\alpha$ and $\gamma$, if $\alpha<\gamma$, then
$\mIter\alpha \hple \mIter\gamma$.
\end{prop}
\begin{proof}
The proposition is proved by induction on $\gamma$ for all $\alpha<\gamma$.
Note that it is true when $\alpha=0$, even
though $\mIter 0$ is not itself hereditarily positive.
If $\gamma = 0$, then the claim is vacuously true.

Suppose that $\gamma = \delta+1$ and fix any $\alpha<\gamma$.  
By the induction hypothesis
$\mIter\alpha\hple \mIter\delta$, so it suffices to show that
$\mIter\delta\hple\mIter\gamma$.  To do so, fix h.p.\ functionals $f$ and
$x$.  Since
$f$ and $\mIter\delta fx$ are h.p. (notice that this is true even when
$\delta = 0$, since then $\mIter\delta fx = x$), $\mIter\delta fx\hple f(\mIter\delta fx) = \mIter\gamma fx$.
Since $f$ and $x$ were chosen arbitrarily,
$\mIter\delta\hple\mIter\gamma$.

Suppose that $\gamma$ is a limit and fix any $\alpha<\gamma$.  
By the induction hypothesis the
sequence $\setidx{\mIter\xi}{\xi<\gamma}$ is non-decreasing with
respect to $\hple$.  Thus, by 
Lemma~\ref{lem:hple-facts}(\ref{item:limsup-eq-sup}),
$\mIter\gamma = \limsup_{\xi\to\gamma}\set{\mIter\xi}\hpeq
\sup_{\xi<\gamma}\set{\mIter\xi}$.
Since $\alpha<\gamma$, there is some
$\delta<\gamma$ such that $\alpha<\delta$, which, by the induction hypothesis
applied to $\delta$, implies that
$\mIter\alpha\hple\mIter\delta\hple\sup_{\xi<\gamma}\set{\mIter\xi}
\hpeq\mIter\gamma$.
\end{proof}

At this point, we are almost done, because if $\setidx{\alpha_\xi}{\xi<\zeta}$
is an increasing sequence of ordinals, then
$\limsup_{\xi\to\zeta}\set{\mIter{\alpha_\xi}} = 
\sup_{\xi<\zeta}\set{\mIter{\alpha_\xi}}$ (recall that the difficulty was
evaluating the $\limsup$ over a subsequence).
But first we need to ensure that the supremum
is itself an iteration functional.  With a little extra effort,
we can prove a more general result:
$\limsup_{\xi\to\zeta}\set{\mIter{\alpha_\xi}}$ is an iteration 
functional for \emph{any} sequence of 
ordinals~$\setidx{\alpha_\xi}{\xi<\zeta}$.  To prove this, we
combine Prop.~\ref{prop:miter-incr} with the fact that the $\limsup$ of
a sequence of ordinals can always be calculated as the supremum over some
tail of the sequence.

\begin{lemma}
\label{clm:limsup-like-sup}
For any sequence of ordinals
$\setidx{\alpha_\xi}{\xi<\zeta}$, there is an ordinal $\mu<\zeta$ such that
$\limsup_{\xi\to\zeta}\{\alpha_\xi\} = \sup_{\mu\leq\xi<\zeta}\{\alpha_\xi\}$.
\end{lemma}
\begin{proof}
By definition, $\limsup_{\xi\to\zeta}\set{\alpha_\xi} = 
\mathop{\vphantom{\sup}\inf}_{\gamma<\zeta}\set{\sup_{\gamma\leq\xi<\zeta}\set{\alpha_\xi}}$.
Since
any set of ordinals attains its infimum, there is some $\mu<\zeta$ such that
$\mathop{\vphantom{\sup}\inf}_{\gamma<\zeta}\set{\sup_{\gamma\leq\xi<\zeta}\set{\alpha_\xi}} =
\sup_{\mu\leq\xi <\zeta}\set{\alpha_\xi}$.
\end{proof}

An analogous fact holds for sequences of iteration functionals:

\begin{lemma}
\label{clm:miter-limsup-sup}
For any sequence of ordinals $\{\alpha_\xi\}_{\xi<\zeta}$, 
take $\mu$ as in Lemma~\ref{clm:limsup-like-sup}; then
$\limsup_{\xi\to\zeta}\{\mIter{\alpha_\xi}\} \hpeq
\sup_{\mu\leq\xi<\zeta}\{\mIter{\alpha_\xi}\}$.
\end{lemma}
\begin{proof}
By the  choice of $\mu$, we have
$\limsup_{\xi\to\zeta}\{\mIter{\alpha_\xi}\} =
\mathop{\vphantom{\sup}\inf}_{\gamma<\zeta}\bigl\{\sup_{\gamma\leq\xi<\zeta}\{\mIter{\alpha_\xi}\}\bigr\}
\leq \sup_{\mu\leq\xi<\zeta}\{\mIter{\alpha_\xi}\}$.
For the reverse inequality, fix any $\delta$ such that $\mu\leq\delta<\zeta$;
then $\alpha_\delta \leq
\sup_{\mu\leq\xi<\zeta}\{\alpha_\xi\} = 
\inf_{\gamma<\zeta}\bigl\{\sup_{\gamma\leq\xi<\zeta}\{\alpha_\xi\}\bigr\}$
by the choice of $\mu$.
So for any 
$\gamma<\zeta$, $\alpha_\delta\leq\sup_{\gamma\leq\xi<\zeta}\set{\alpha_\xi}$,
and therefore there is some $\xi_\gamma\geq\gamma$ such that $\alpha_\delta
\leq\alpha_{\xi_\gamma}$, which by
Prop.~\ref{prop:miter-incr} implies that
$\mIter{\alpha_\delta} \hple \mIter{\alpha_{\xi_\gamma}}$.  
Keeping in mind that $\delta$ is fixed while $\gamma$ was chosen
arbitrarily,
$\mIter{\alpha_\delta}\hple
\limsup_{\gamma\to\zeta}\{\mIter{\alpha_{\xi_\gamma}}\} \hple
\limsup_{\xi\to\zeta}\{\mIter{\alpha_\xi}\}$; the final
inequality follows from
Lemma~\ref{lem:hple-facts}(\ref{item:boundedness}).
Since $\delta$ was chosen arbitrarily between $\mu$ and $\zeta$, 
this implies that
$\sup_{\mu\leq\xi<\zeta}\{\mIter{\alpha_\xi}\}\hple
\limsup_{\xi\to\zeta}\{\mIter{\alpha_\xi}\}$.
\end{proof}

\begin{prop}
\label{clm:miter-subseq}
For any sequence of ordinals $\setidx{\alpha_\xi}{\xi<\zeta}$,
$\limsup_{\xi\to\zeta}\{\mIter{\alpha_\xi}\}\hpeq 
\mIter{\limsup\set{\alpha_\xi}}$.
\end{prop}
\begin{proof}
Fix $\mu$ as in Lemma~\ref{clm:limsup-like-sup} and set
$\alpha\eqdef\limsup\set{\alpha_\xi} = 
\sup_{\mu\leq\xi<\zeta}\set{\alpha_\xi}$.
First, suppose that for all $\gamma$ there is $\xi_\gamma\geq\gamma$ such that
$\alpha_{\xi_\gamma} = \alpha$.  Then since $\alpha_\xi\leq\alpha$ for
all $\mu\leq\xi<\zeta$, $\sup_{\mu\leq\xi<\zeta}
\set{\mIter{\alpha_\xi}}\hple\mIter\alpha$.  
On the other hand, $\mu\leq\xi_\mu<\zeta$ and
$\alpha_{\xi_\mu} = \alpha$, so 
$\mIter\alpha\hple\sup_{\mu\leq\xi<\zeta}\set{\mIter{\alpha_\xi}}$, and
therefore $\limsup\set{\mIter{\alpha_\xi}}\hpeq
\sup_{\mu\leq\xi<\zeta}\set{\mIter{\alpha_\xi}} \hpeq\mIter\alpha$.

If there is some $\gamma$ such that $\alpha_\xi<\alpha$ for all
$\xi\geq\gamma$, then we can still conclude that
$\limsup_{\xi\to\zeta}\{\mIter{\alpha_\xi}\} \hple\mIter\alpha$.
Note that in this situation, $\alpha$ must be a limit.  To show that
the reverse inequality holds, fix any $\gamma<\alpha$; then there
is $\delta\geq \mu$ such that $\gamma\leq\alpha_\delta$, so by
Prop.~\ref{prop:miter-incr}
$\mIter\gamma\hple\mIter{\alpha_\delta}\hple
\sup_{\mu\leq\xi<\zeta}\{\mIter{\alpha_\xi}\} \hpeq
\limsup_{\xi\to\zeta}\{\mIter{\alpha_\xi}\}$.
But since $\gamma$ was chosen arbitrarily, this implies that
$\mIter\alpha =
\limsup_{\gamma\to\alpha}\{\mIter\gamma\} \hple
\limsup_{\xi\to\zeta}\{\mIter{\alpha_\xi}\}$.
The inequality follows from
Lemma~\ref{lem:hple-facts}(\ref{item:boundedness}) by considering
$\limsup\set{\mIter{\alpha_\xi}}$ as a single h.p.\ functional
bounding \emph{all} of the $\mIter\gamma$.
\end{proof}

Now we
arrive at the main result relating compositions of functionals of the form
$\mIter\alpha$ to ordinal arithmetic:

\begin{thm}[Iteration Functionals under $\hpeq$]
\label{thm:comp-of-miter}
Let $f\in\Omega_{\tau\arrow\tau}$ and $x\in\Omega_\tau$
be h.p.  Then for any countable $\alpha$ and $\gamma$:
\begin{enumerate}
\item \label{item:add}
$\mIter\alpha^\tau f(\mIter\gamma^\tau fx) \hpeq 
\mIter{\gamma+\alpha}^\tau fx$.
\item \label{item:mult}
$\mIter\alpha^\tau (\mIter\gamma^\tau f) \hpeq \mIter{\gamma\alpha}^\tau f$.
\item \label{item:exp}
$\mIter\alpha^{\tau\arrow\tau}(\mIter\gamma^{\tau})\hpeq
      \mIter{\gamma^\alpha}^\tau$.
\end{enumerate}
\end{thm}
\begin{proof}
Each part is proved
by induction on $\alpha$; we do (\ref{item:mult}) as an example.
If $\alpha = 0$ and $x$ is h.p., then $\mIter\alpha(\mIter\gamma f)x =
x = \mIter{\gamma\alpha}fx$.  

If $\alpha = \delta+1$, then
$\mIter\alpha(\mIter\gamma f)x = 
  (\mIter\gamma f)\bigl(\mIter\delta(\mIter\gamma f)x\bigr) \hpeq
  (\mIter\gamma f)\bigl(\mIter{\gamma\delta}fx\bigr) \hpeq
  (\mIter{\gamma\delta+\gamma}fx) \hpeq
  (\mIter{\gamma(\delta+1)}fx)
$.
The second equality is the induction hypothesis and the third is an
application of (\ref{item:add}).

If $\alpha$ is a limit, then
$\mIter\alpha(\mIter\gamma f) \hpeq
  \limsup_{\xi\to\alpha}\{\mIter\xi(\mIter\gamma f)\} \hpeq
  \limsup_{\xi\to\alpha}\{\mIter{\gamma\xi}f\}\hpeq
  \mIter{\gamma\alpha} f
$,
with the middle equality following from the induction hypothesis and
the last one by Prop.~\ref{clm:miter-subseq}.
\end{proof}

It is useful to note why $f$ and $x$ are required to be h.p.\ in 
Theorem~\ref{thm:comp-of-miter}.  In the last equality of the limit
case, we use Prop.~\ref{clm:miter-subseq}, which asserts only that
$\limsup_{\xi\to\alpha}\set{\mIter{\gamma\xi}}\hpeq\mIter{\gamma\alpha}$.
Thus, \emph{when $f$ is h.p.}, we can conclude that
$$
  \limsup_{\xi\to\alpha}\set{\mIter{\gamma\xi}f} =
  \Bigl(\limsup_{\xi\to\alpha}\set{\mIter{\gamma\xi}}\Bigr)(f)\hpeq
  \mIter{\gamma\alpha}f
$$
In particular, the ``alternating'' function which we considered
in Section~\ref{sec:Preliminaries} is not itself h.p., and this last
argument would fail for that function.

\bibliographystyle{plain}
\bibliography{general,lambda-bib,misc}

\end{document}